\documentclass[]{article}

\usepackage{graphicx}
\usepackage{hyperref}
\usepackage[utf8]{inputenc}

\title{                                     
Comprehensive passenger demand-dependent traffic control on a metro
line with a junction and a derivation of the traffic phases}

\author{
  \textbf{Florian Schanzenbacher}\\
  RATP, Enpc and Ifsttar\\
  40 bis, rue Roger Salengro, 94724 Fontenay-sous-Bois cedex France\\
  Tel: +331 58 76 45 50; Email: florian.schanzenbacher@ratp.fr\\
  ORCiD number: 0000-0002-0876-1715\\
  \hfill\break\\% this is a way to add line numbering on empty line
  \textbf{Nadir Farhi}\\
  Ifsttar / Cosys / Grettia\\
  14-20 Boulevard Newton, 77447 Marne-la-Vall\'ee cedex 2 France\\
  Tel: +331 81 66 87 04; Email: nadir.farhi@ifsttar.fr\\
  \hfill\break\\% this is a way to add line numbering on empty line
  \textbf{Fabien Leurent}\\
  Enpc / Lvmt\\
  6-8 Avenue Blaise Pascal, 77455 Marne-la-Vall\'ee cedex 2 France\\
  Tel: +331 81 66 88 54; Email: fabien.leurent@enpc.fr\\
  \hfill\break\\% this is a way to add line numbering on empty line
  \textbf{G\'erard Gabriel}\\
  RATP / MOP / GEF\\
  40 bis, rue Roger Salengro, 94724 Fontenay-sous-Bois cedex France\\
  Tel: +331 58 76 48 57; Email: gerard.gabriel@ratp.fr
}

\begin{document}
\maketitle

\section{Introduction and review of preceding works}
Mass transit metro lines of world cities like Paris have to cope with a high passenger travel demand.
A recurrent phenomenon on these mass transit systems is the amplification of small perturbations on the train time-headway.
Due to the complexity of the system, minor perturbations occur frequently, resulting in deviations on the train headway.
In case of a high demand, even small disturbances may cause a serious accumulation of passengers on the platforms,
which results in longer dwell times and provokes a cascade effect, i.e. the perturbation quickly propagates through the line.
We present here a discrete event traffic model of the train dynamics on a metro line with a junction, taking into account the passenger travel demand.
This model allows a complete understanding of the traffic phases of the train dynamics, incorporating passenger travel demand, i.e. it allows to study the above described problem.

The following works are inspired by some pioneer works on traffic control for metro loop lines by~\cite{Bres91,Fer05}.
The authors of~\cite{Bres91} have presented a model for a control of metro loop line train dynamics while the authors of~\cite{Fer05} have proposed a model predictive control approach to minimize train time-headway variance. The model predictive control approach has been applied by~\cite{SCF16} to a dense rail commuter line in the area of Paris.

We distinguish three parts of the metro line: a central part and two branches crossing at the junction.
The metro line is discretized in space into a number of segments as in~\cite{FNHL16,FNHL17a, FNHL17b, F18, SFCLG17, SFLG18a, SFLG18b}, all following the same discrete event modeling approach.
A first model for the train dynamics on a metro line with a junction has been proposed in~\cite{SFCLG17}, with train dwell and run times respecting given lower bounds, i.e. they are \textit{constant} and \textit{independent of the passenger travel demand}.
A second model has recently been proposed in~\cite{SFLG18a},
with passenger demand-dependent dwell and run times, but in a metro \textit{loop} line \textit{without junction}.
In both cases, it has been shown that the train dynamics
admit a stationary regime with a unique asymptotic average train time-headway.
Moreover, the latter is derived analytically as a function of the number of running trains,
and, for the line \textit{with a junction}, of the difference between the number of running trains on each of the two branches.

The model we present in this article extends both models of~\cite{SFCLG17} and~\cite{SFLG18a} in the following way.
We show that the results on the stability of the train dynamics for a metro line \textit{with a junction} still hold under the following changes.
First, we consider train dwell times as a function of the passenger travel demand and the train headway, bounded by a maximum dwell time.
Second, we introduce a run time control which compensates eventual extensions of the dwell time.
This permits to guarantee train dynamics stability.
We derive analytically the phase diagrams of the train dynamics
on a metro line with a junction and with a demand-depend traffic control.
The analytic formulas and diagrams can be used to control the train passing order at the junction, in the case
where the system switches between two steady states, for example, due to a change in passenger demand.
The derivation of the traffic phases of the train dynamics in all the models is based on a Max-plus algebra modeling approach.
For more details on theoretical background of Max-plus algebra systems and the main theorem underlying the derivation of the dynamics, the reader is referred to~\cite{BSOQ92,CCGMQ98}.

Finally,~\cite{SFLG18b} propose a new control on the train dwell times, which allows to minimize headway variance. However, this control is not taken into account in the here presented model for a metro line \textit{with a junction}.

% =======================================================
\section{Review of the modeling approach}
% =======================================================
\subsection{Hypotheses}
\begin{figure}[!ht]
  \centering
  \includegraphics[width=0.6\textwidth]{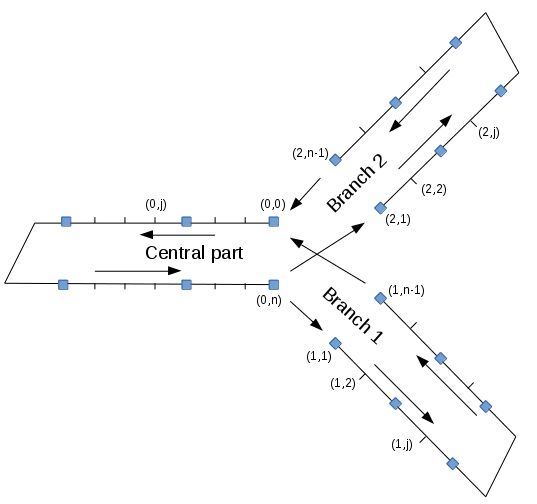}
  \caption{Schema of a metro line with one junction and the here proposed notation.}\label{fig_line}
\end{figure}

We consider a metro line with a symmetrically operated junction as in~\cite{SFCLG17}.
The line is discretized into segments and our main variables are the $k^{th}$ departures
from each segment $j$ on the three parts $u$ of the line, see FIGURE~\ref{fig_line}.
We note that there can be more than one segment between two platforms and that there can be maximum one train in a segment at a time.
Trains are running between one of the terminus on the branches and the one on the central part and stop at every platform.
There are no further constraints on the layout of the line, especially the branches can have different lengths and different number of platforms.

Passenger arrival rates to the platforms are supposed to be constant over a time interval which can be chosen in an appropriate way.
Passenger boarding and alighting rates are constant, i.e. there is no congestion effect in case of crowded vehicles and platforms.
Finally, we consider full observability and controllability of the system.

\subsection{Modeling}
Our model describes the train dynamics by taking into account the traffic constraints related to train running, dwell and safe separation times.
Moreover, the model includes two control laws on the train dwell and run times, as functions of the train time headways and of the passenger travel demand levels at every platform.
The dwell time control law proposes train dwell times as long as the time needed for passenger alighting and boarding at every platform, bounded by a maximum dwell time.
The run time control law proposes a compensation of an eventual extension of a train dwell time at a given platform, by
shortening the running time of the same train on the downstream segments.

\section{Findings: Derivation of the traffic phases of the demand-dependent train dynamics in a metro line with a junction}

The results presented in the following, are based on the main result of~\cite{SFCLG17} and adapted to the case where the train dwell times depend on the
passenger travel demand, and where the running times are controlled.
We consider the following notations.
\begin{table}[!ht]
	\caption{Parameters of the traffic phases of the train dynamics}%\label{tab:notations_2}
	\begin{center}
		\begin{tabular}{l l}
		\hline
		$m$ & the total number of trains on the line.\\
  		$\Delta m$ & $= m_2 - m_1$ the difference in the number of trains between branches 2 and 1.\\ \hline
		\end{tabular}
	\end{center}
\end{table}

\begin{figure}[!ht]
  \centering
  \includegraphics[width=0.6\textwidth]{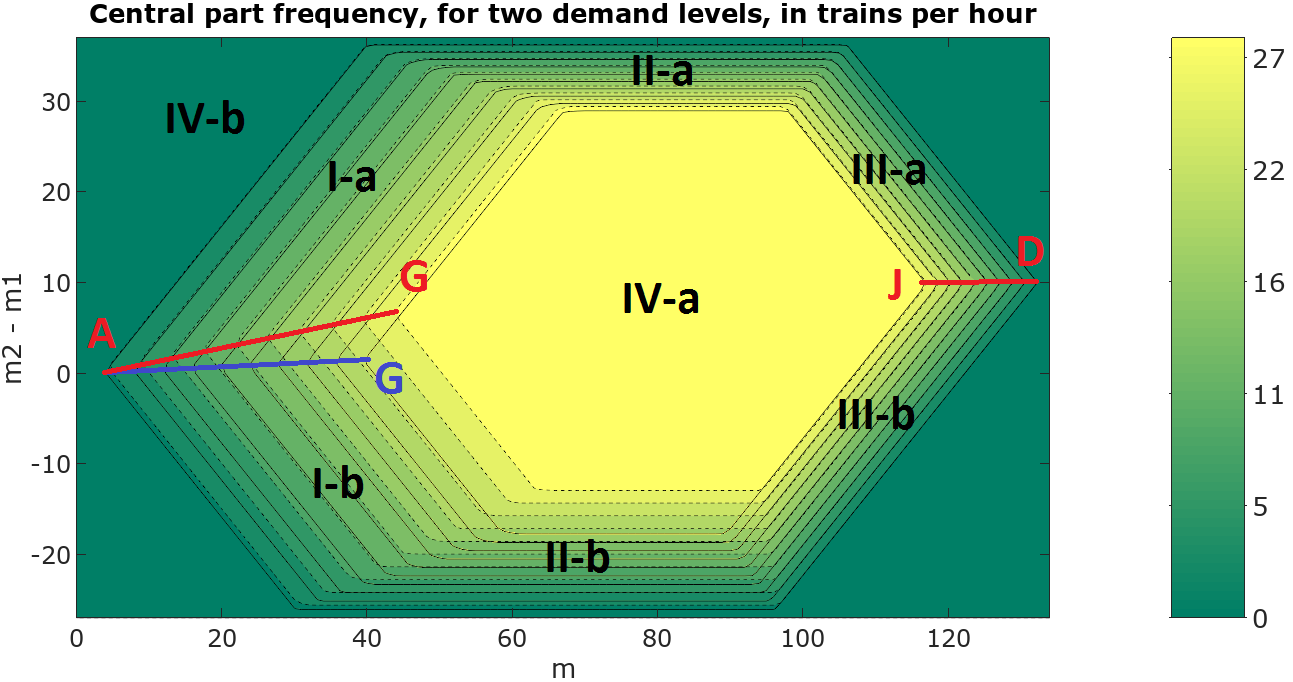}\label{fig_surface}
  \caption{Average central part frequency $f_0$ on RATP metro line 13, Paris. The line consists of two branches and a central part, crossing at a junction.
  Depicted are the eight traffic phases of the train dynamics.
  First, today's peak hour demand (solid lines, blue line [AG]).
  Second, a demand increase going exclusively on branch 2 and on the central part (dashed lines, red line [AG]).
  }
\end{figure}
FIGURE~2 depicts the main result of the discrete event traffic model of a metro line with a junction and demand-dependent dwell times and a run time control, applied to metro line 13 of Paris.
There are eight traffic phases of the train dynamics, given here by the average frequency on the central part of the line, which is twice the one on the branches. It is depicted on the color scale and has its maximum \textit{on Paris metro line 13} at around 27 train per hour.
The average frequency on the line with a junction depends on two parameters:
The total number of trains $m$,
the difference between the number of trains between the two branches, $m_2 - m_1$.
There are eight traffic phases in total, i.e. two free flow phases (I-a/I-b), two unbalanced branches phases (II-a/II-b),
two congested traffic phases (III-a/III-b) as well as a zero and maximum frequency phase (IV-a/IV-b).
Finally, the straight lines [AG] and [JD] demonstrate the influence of the passenger demand.
Whereas [JD], separating the two congested traffic phases, is independent of the passenger demand,
[AG], separating the two free flow traffic phases, deviates with a changing passenger demand deviation between the two branches.
The blue line and the solid lines represent today's peak hour demand on line 13.
The red line and the dashed lines correspond to an example where demand increases exclusively on branch 2 and on the central part, i.e. demand on branch 1 remains unchanged.

\section{Conclusion}
We have presented a discrete event traffic model for a metro line with a symmetrically operated junction.
The model describes the train dynamics by means of time constraints on the train run, dwell and safe separation 
times. The train dwell times are controlled in function of the passenger demand at the platforms.
The train run times on inter-stations are controlled in such a way that they compensate eventual extensions
of the train dwell times at the upstream platforms, which may due to high levels of the passenger demand.
Under sufficient margins on the train run times, the train dynamics attain
a stable stationary regime.
The asymptotic average train time-headway is given as
a function of the total number of trains, of the difference between the number of trains
on the two branches, and of the level of passenger travel demand.
We thus obtained the traffic phases of the train dynamics.
The results have been presented in a descriptive 2D graphic which illustrates the influence
of the demand on the train dynamics of the metro line.

\end{document}